\pdfoutput=1
\RequirePackage{ifpdf}
\ifpdf 
\documentclass[pdftex]{sigma}
\else
\documentclass{sigma}
\fi

\numberwithin{equation}{section}
\numberwithin{theorem}{section}

\let\frak\mathfrak
\let\Bbb\mathbb

\usepackage{eucal}

 \let\alb\allowbreak
\def\fratop{\genfrac{}{}{0pt}1}

\newtheorem{thm}{Theorem}[section]

\newtheorem{lem}[thm]{Lemma}

\theoremstyle{definition}
\newtheorem*{rem}{Remark}
\newtheorem*{example}{Example}

 \let\eps\varepsilon \let\epsilon\eps

\let\la\lambda 

\let\si\sigma 
 
 \let\phi\varphi

\let\Tilde\widetilde

\let\der\partial

\let\geq\geqslant
\let\le\leqslant
\let\leq\leqslant

\def\C{\Bbb C}

\def\Z{\Bbb Z}

\def\gl{\frak{gl}}

\def\lsym#1{#1\alb\dots\relax#1\alb} \def\lc{\lsym,}

\let\on\operatorname

\def\diag{\on{diag}}

\def\rdet{\on{rdet}}
\def\Res{\mathop{\mathrm{Res}}\limits}

\def\Wr{\on{Wr}}

\def\gln{\gl_N}

\let\bs\boldsymbol

\let\nc\newcommand

\nc{\pone}{\C{\mathbb P}^1}

\def\0{{\>0}}
\def\ss{{\<\>\raise.22ex\hbox{$\sssty S$}}}
\def\yy{{\<\>\raise.22ex\hbox{$\sssty Y\!\!\!$}}}

\def\sing{\on{Sing\>}\nolimits}

\def\%#10{#1,\<\>0}
\def\^#1{^{[#1]}}
\def\@#1{^{(#1)}}
\def\+#1{^{\{#1\}}}

\def\KZ/{{\sl KZ\/}}
\def\qKZ/{{\sl qKZ\/}}

\let\co c
\let\Co C
\let\Hh K

\let\rr r
\let\yh h
\nc{\bla}{{\bs\la}}
\nc{\D}{{\mathcal D}}
\nc{\W}{{\mathcal W}}
\nc{\Pa}{{\mathcal P}}
\let\yp p
\def\glN{\gl_N}
\let\kk K
\nc{\TT}{{\bs t}}
\nc{\zz}{{\bs z}}
\nc{\Spect}{\on{Spec}\nolimits}
\nc{\Spec}{\Spect_\bla}
\nc{\Crit}{\on{Crit}}

\begin{document}

\renewcommand{\PaperNumber}{072}

\FirstPageHeading

\ShortArticleName{KZ Characteristic Variety as the Zero Set of Classical Calogero--Moser
Hamiltonians}

\ArticleName{KZ Characteristic Variety as the Zero Set\\ of Classical
Calogero--Moser Hamiltonians}

\AuthorNameForHeading{E.~Mukhin, V.~Tarasov, and A.~Varchenko}

\Author{Evgeny MUKHIN~$^\dag$, Vitaly~TARASOV~$^{\dag\ddag}$ and Alexander~VARCHENKO~$^\S$}

\Address{$^\dag$~Department of Mathematical Sciences,
Indiana University~-- Purdue University Indianapolis,\\
\hphantom{$^\dag$}~402 North Blackford St, Indianapolis, IN 46202-3216, USA}
\EmailD{\href{mailto:mukhin@math.iupui.edu}{mukhin@math.iupui.edu}, \href{mailto:vtarasov@math.iupui.edu}{vtarasov@math.iupui.edu}}

\Address{$^\ddag$~St.~Petersburg Branch of Steklov Mathematical Institute,\\
\hphantom{$^\ddag$}~Fontanka 27, St.~Petersburg, 191023, Russia}
\EmailD{\href{mailto:vt@pdmi.ras.ru}{vt@pdmi.ras.ru}}

\Address{$^\S$~Department of Mathematics, University of North Carolina
at Chapel Hill,\\
\hphantom{$^\S$}~Chapel Hill, NC 27599-3250, USA}
\EmailD{\href{mailto:anv@email.unc.edu}{anv@email.unc.edu}}

\ArticleDates{Received June 09, 2012, in f\/inal form October 04, 2012; Published online October 16, 2012}

\Abstract{We discuss a relation between the characteristic variety of the KZ equations
and the zero set of the classical Calogero--Moser Hamiltonians.}


\Classification{82B23; 17B80}

\section{Statement of results}
\label{intro}

\subsection{Motivation}

This paper is motivated by the Givental and Kim observation~\cite{GK}
that the characteristic variety of the quantum dif\/ferential equation
of a f\/lag variety is a Lagrangian variety of the classical Toda lattice.
The quantum dif\/ferential equation is a system of dif\/ferential equations
$\hbar \der_i\psi = b_i\circ \psi$, $i=1,\dots,r,$
def\/ined by the quantum multiplication $\circ$ and
depending on a parameter~$\hbar$. The system def\/ines a f\/lat connection for all nonzero
values of~$\hbar$. Givental and Kim, in particular, observe that the characteristic variety
of this system is the Lagrangian variety of the classical Toda lattice, def\/ined
by equating to zero the f\/irst integrals of the Toda lattice.

In this paper we describe a similar relation between the KZ equation and the classical Calogero--Moser system.
On numerous relations between the KZ equations and quantum Ca\-lo\-ge\-ro--Moser systems see \cite{C1,C2,FV2,FV1,M}.

\subsection{Classical Calogero--Moser system}
\label{sec Calogero--Moser system}

Fix an integer $n\geq 2$. Denote
$\Delta =\{
\bs z=(z_1,\dots,z_n) \in \C^n\ | \ z_a=z_b\ \text{for some}\ a\neq b\}$, the union of diagonals.
Consider the cotangent bundle $T^*(\C^n-\Delta)$ with symplectic form
$\omega = \sum\limits_{a=1}^n dp_a\wedge dz_a$, where $p_1,\dots,p_n$ are coordinates on f\/ibers.
The classical Calogero--Moser system on $T^*(\C^n-\Delta)$ is def\/ined by the Hamiltonian
\begin{gather*}
\mathcal H   =  \sum_{a=1}^n  p_a^2 -  \sum_{1\leq a<b\leq n } \frac 2{(z_a-z_b)^2}.
\end{gather*}
The system is completely integrable.
For
\begin{gather}
\label{Q}
Q =
\begin{pmatrix}
p_1 & \dfrac{1}{z_1-z_2} & \dfrac{1}{z_1-z_3} &  \dots & \dfrac{1}{z_1-z_n}
\\[5pt]
\dfrac{1}{z_2-z_1} & p_2 & \dfrac{1}{z_2-z_3} &  \dots & \dfrac{1}{z_2-z_n}
\\[5pt]
{}\dots & {}\dots & {}\dots & \,\dots & {}\dots
\\[5pt]
\dfrac{1}{z_n-z_1} & \dfrac{1}{z_n-z_2}& \dfrac{1}{z_n-z_3}& \,\dots & p_n
\end{pmatrix} ,
\end{gather}
let $\det(u-Q)= u^n - Q_1u^{n-1}+\cdots \pm Q_n$ be the characteristic polynomial.
Then $Q_1,\dots,Q_n$ is a complete list of commuting f\/irst integrals,
and
$\mathcal H= Q_1^2-Q_2$.

We will be interested in the subvariety $L_0\subset T^*(\C^n-\Delta)$ def\/ined
by the equations
\begin{gather}
\label{L0}
L_0 = \big\{(\bs z,\bs p)\in T^*(\C^n-\Delta)\; |\; Q_a(\bs z, \bs p)=0, \
a=1,\dots,n\big\} .
\end{gather}

\begin{thm}[\cite{Wi}]
\label{thm L_0 smooth}
For any $n$, the subvariety $L_0$ is smooth and Lagrangian.
\end{thm}

See propositions in Section~6 of \cite{Wi}. Another proof of
Theorem \ref{thm L_0 smooth} will be given in Section~\ref{sec part ii proof}.

\subsection{Gaudin Hamiltonians and KZ characteristic variety}
\label{Gaudin Hamiltonians}

Fix an integer $N\geq 2$.
Denote $V=\C^N$ the vector representation of the Lie algebra $\glN$.
The Hamiltonians of the quantum Gaudin model are the linear operators
$H_1\lc H_n$ on the space~$V^{\otimes n}$,
\begin{gather}
\label{Ha}
H_a(\bs z)   =
\sum_{i,j=1}^N \sum_{b\neq a}  \frac{e_{ij}^{(a)}e_{ji}^{(b)}}{z_a-z_b} ,
\end{gather}
where  $e_{ij}$ are the standard generators of $\gln$,  $e_{ij}^{(a)}$
is the image of $1^{\otimes(a-1)}\otimes e_{ij}\otimes1^{\otimes(n-a)}$, and
$z_1,\dots,z_n$ are distinct complex numbers, see~\cite{G}. The operators commute,
$[H_a(\bs z),H_b(\bs z)]=0$ for all~$a$,~$b$.
The operators commute with the $\gln$-action on~$V^{\otimes n}$.

Let $\bla = (\la_1,\dots,\la_N) \in \Z^N_{\geq 0}$ be a partition of $n$ with at most $N$ parts,
$\la_1\geq\dots\geq\la_N$, $|\bla|=\la_1+\dots+\la_N=n$. Denote
\begin{gather*}
\sing V^{\otimes n}[\bla] = \big\{v\in V^{\otimes n}\; | \; e_{ii}v=\la_iv, \  i=1,\dots,N; \   e_{ij}v=0\ \text{for all}\ i<j\big\},
\end{gather*}
the subspace of singular vectors of weight $\bla$. The Gaudin Hamiltonians preserve
$\sing V^{\otimes n}[\bla]$.
We def\/ine the spectral variety of the Gaudin model on $\sing V^{\otimes n}[\bla]$,
\begin{gather*}
\Spect_{N,\bla} = \big\{(\bs z,\bs p)\in T^*(\C^n-\Delta)\; | \; \exists\,
v\in \sing\,V^{\otimes n}[\bla]\ \text{with}\ H_a(\bs z)v = p_av,\ a=1,\dots,n\big\}.
\end{gather*}
The spectral variety is a Lagrangian subvariety of $T^*(\C^n-\Delta)$, see,
for example, Proposition~1.5 in \cite{RV}.

The Gaudin Hamiltonians are the right hand sides of the KZ equations,
\[
\kappa \der_{z_i} I(\bs z) = H_i(\bs z) I(\bs z), \qquad i=1,\dots,n,
\]
where $\kappa\in\C^\times$ is a parameter.
The spectral variety $\Spect_{N,\bla}$ is, by def\/inition, the characteristic
variety of the $\kappa$-dependent $D$-module def\/ined by the KZ equations with
values $\sing V^{\otimes n}[\bla]$.

\begin{example}
If $\bla=(n,0,\dots,0)$, then $\Spect_{N,\bla}$ is given by the equations
$p_a = \sum\limits_{b\neq a}(z_a-z_b)^{-1}$,  $a=1\lc n$.
If $N=n$ and $\bla=(1,\dots,1)$, then $\Spect_{N,\bla}$ is given by the equations
$p_a = -\sum\limits_{b\neq a}(z_a-z_b)^{-1}$,  $a=1\lc n$.
\end{example}

\begin{thm}\label{thm spec-lambda in L_0}\qquad\null
\begin{enumerate}\itemsep=0pt
\item[$(i)$]
The variety $\Spect_{N,\bla}$ does not depend on $N$. Namely, consider $\bla$ as a partition of $n$ with
at most $N+1$ parts. Then $\Spect_{N,\bla}=\Spect_{N+1,\bla}$.
From now on we denote $\Spect_{N,\bla}$ by $\Spec$.
\item[$(ii)$] For any~$n$, the variety $L_0$ is the disjoint union of the
varieties~$\Spec$, where the union is over all partitions~$\bla$ of~$n$.
\end{enumerate}
\end{thm}

By Theorem~\ref{thm L_0 smooth}, each $\Spec$ is smooth and Lagrangian.

Part $(i)$ of Theorem \ref{thm spec-lambda in L_0} is proved in Section \ref{Proof of part i}, and part $(ii)$ is proved in
Section~\ref{sec part ii proof}.

\subsection[Master function generates ${\rm Spec}_{\lambda}$]{Master function generates $\boldsymbol{{\rm Spec}_{\lambda}}$}
\label{Generating function of Spec-bla}

Let $\bs\la$ be a partition of $n$ with at most $N$ parts.
Denote  $l_a=\sum\limits_{b=a+1}^N \la_b$, $a=1,\dots,N-1$.
Denote $l=l_1+\dots +l_{N-1}$. Consider the set of $l$ variables
\[
\bs t =\big(t^{(1)}_{1},\dots,t_{l_1}^{(1)},\dots,
t^{(N-1)}_{1},\dots,t_{l_{N-1}}^{(N-1)}\big)
\]
and the af\/f\/ine space $\C^n\times\C^l$ with coordinates
$\bs z$, $\bs t$. The function $\Phi_{N,\bla}: \C^n\times\C^l\to\C$,
\begin{gather*}
\Phi_{N,\bla} (\bs z,\bs t) =
\sum_{1\leq a<b\leq n} \log(z_a-z_b)-\sum_{a=1}^{n}\sum_{i=1}^{l_a}
\log\big(t_i^{(1)}-z_a\big)
\\
\hphantom{\Phi_{N,\bla} (\bs z,\bs t) = }{} +
2 \sum_{k=1}^{N-1} \sum_{1\leq i<j\leq l_k} \log\big(t_i^{(k)}-t_j^{(k)}\big)-
\sum_{k=0}^{N-2} \sum_{i=1}^{l_k} \sum_{j=1}^{l_{k+1}}
\log\big(t_i^{(k)}-t_j^{(k+1)}\big)
\end{gather*}
is called the master function, see \cite{SV2, SV1}.

The master function depends on $\bla$, but not on $N$.
Namely, consider $\bla$ as a partition of $n$ with at most $N+1$ parts.
Then $\Phi_{N,\bla}=\Phi_{N+1,\bla}$. From now on we denote $\Phi_{N,\bla}$
by~$\Phi_\bla$.

Critical points of $\Phi_\bla$ with respect to $\bs t$ are given by the
equation $d_{\bs t}\Phi_\bla = 0$. Denote by~$\Crit_\bla$ the critical set of
$\Phi_\bla$ with respect to~$\bs t$,
\[
\Crit_\bla = \big\{(\bs z,\bs t)\in \C^n\times\C^l\; | \; d_\TT\Phi_\bla(\bs z,\bs t) = 0\big\} .
\]
This is an algebraic subset of the domain of $ \C^n\times\C^l$, where the master function
is a regular (multivalued) function.
Denote by $L_\bla\subset T^*(\C^n-\Delta)$ the image of the map
\[
\Crit_\bla \to T^*(\C^n-\Delta), \qquad (\bs z,\bs t) \mapsto (\bs z,\bs p), \qquad
\text{where}
\quad p_a = \frac {\der \Phi_\bla}{\der z_a}(\bs z,\bs t),\quad a=1,\dots,n.
\]

\begin{thm}
\label{thm gener fun}
For any $n$ and a partition $\bla$ of $n$, we have $L_\bla\subset\Spec$ and
$\Spec$ is the closure of~$L_\bla$ in $T^*(\C^n-\Delta)$.
\end{thm}

Theorem \ref{thm gener fun} is proved in Section \ref{Proof of Theorem gener fun}.

\subsection[Calogero-Moser space $C_n$ and cotangent bundle $T^*(\C^n-\Delta)$]{Calogero--Moser space $\boldsymbol{C_n}$ and cotangent bundle $\boldsymbol{T^*(\C^n-\Delta)}$}
\label{sec cot bundle Cal-Moser}

The Calogero--Moser system has singularities if some of $z_1, \dots, z_n$
coincide. These singularities can be resolved and the Calogero--Moser system can
be lifted by the
map $\xi$ given by~\eqref{xi} below
to a regular completely integrable
Hamiltonian system on the Calogero--Moser space $C_n$, see~\cite{KKS}.

Denote
\[
\Tilde C_n = \big\{(Z,Q)\in \frak{gl}_n\times \frak{gl}_n\; |\;
\on{rank}([Z,Q]+1)=1\big\}.
\]
The group $GL_n$ of complex invertible matrices acts on~$\Tilde C_n$
by simultaneous conjugation. The action is free and proper, see~\cite{Wi}.
The quotient space $C_n$ is called the {$n$-th Calogero--Moser space}.
The Calogero--Moser space $C_n$ is a smooth af\/f\/ine variety of
dimension~$2n$, see~\cite{Wi}.

The group $S_n$ freely acts on $\C^n-\Delta$ by permuting coordinates.
The action lifts to a free action on $T^*(\C^n-\Delta)$. Def\/ine the map
\begin{gather}
\label{xi}
\xi : \ T^*\big(\C^n-\Delta\big) \to T^*\big(\C^n-\Delta\big)/S_n \to C_n
\end{gather}
by the rule: $(\zz, \bs p)$ is mapped to $(Z,Q)$,
where $Z=\text{diag}(z_1,\dots,z_n)$ and
$Q$ is def\/ined by~\eqref{Q}.
The map~$\xi$ induces an embedding $T^*(\C^n-\Delta)/S_n \to C_n$
whose image is Zariski open in $C_n$.

Set $\C^{(n)} = \C^n/S_n$ and let
 $\text{spec}(X)\in\C^{(n)}$ stand for the point given by
the eigenvalues of a~square matrix $X$. The canonical map
\begin{gather*}
\pi : \ C_n \to \C^{(n)}\times \C^{(n)}, \qquad (Z,Q) \mapsto (\text{spec}(Z), \text{spec}(Q)),
\end{gather*}
is a f\/inite map of degree $n!$, see~\cite{EG}. This map and its f\/iber over~$0\times 0$ were
studied, for example, in~\cite{EG,FG,Go}.

Let $C_n^0$ be the subvariety $\pi^{-1}(\C^{(n)}\times 0)\subset C_n$.
Identifying $\C^{(n)}\times 0$ with $\C^{(n)}$ we get a map
\begin{gather*}
\pi^0 : \ C^0_n \to \C^{(n)}, \qquad (Z,Q) \mapsto \text{spec}(Z),
\end{gather*}
induced by $\pi$. We will describe $\pi^0$
in Section~\ref{sec description}.

\subsection{Wronski map}

For a partition $\bla$ of $n$, introduce $\bs{\Tilde \lambda} =
\{\Tilde\lambda_1,\dots,\Tilde\lambda_n\}$ by $\Tilde\lambda_i=\la_i+n-i$.
Denote
\begin{gather*}
f_i(u) = u^{\Tilde\lambda_i} +
\sum_{\fratop{j=1}{\Tilde\lambda_i-j\notin \Tilde\bla\!\!}}^{\Tilde\la_i}
f_{ij} u^{\Tilde\lambda_i-j},\qquad i=1,\dots,n.
\end{gather*}
Denote $X_\bla$ the $n$-dimensional af\/f\/ine space of $n$-tuples $\{f_1,\dots,f_n\}$
of such polynomials.
The polynomial algebra
$\C[X_\bla] = \C[f_{ij}, \, i=1,\dots,n,\, j\in\{1,\dots,\Tilde\lambda_i\},\,
\Tilde\lambda_i-j\notin \Tilde\bla]$
is the algebra of regular functions on $X_\bla$.

If $\zz=(z_1,\dots,z_n)$
are coordinates on $\C^n$, then $\bs\si=(\si_1,\dots,\si_n)$,
where $\si_a$ is the $a$-th elementary symmetric function of $z_1,\dots,z_n$,
are coordinates on $\C^{(n)} = \C^n/S_n$.

For arbitrary functions $g_1(u),\dots,g_n(u)$, introduce the Wronskian
determinant by the formula
\[
\Wr(g_1(u),\dots,g_n(u)) =
\det \begin{pmatrix} g_1(u) & g_1'(u) &\dots & g_1^{(n-1)}(u) \\
g_2(u) & g_2'(u) &\dots & g_2^{(n-1)}(u) \\
\dots & \dots &\dots & \dots \\
g_n(u) & g_n'(u) &\dots & g_n^{(n-1)}(u)
\end{pmatrix} .
\]
We have
\begin{gather*}
\Wr(f_1(u),\dots,f_n(u)) =
\prod_{1\le i<j\le n} \big(\Tilde\la_j-\Tilde\la_i\big)\cdot
\left(u^n+\sum_{a=1}^n (-1)^a W_a u^{n-a}\right)
\end{gather*}
with $W_1,\dots,W_n\in\C[X_\bla]$.
Def\/ine an algebra homomorphism
\[
\W_\bla : \ \C[\C^{(n)}] \to \C[X_\bla] ,\qquad
\sigma_a \mapsto W_a .
\]
The corresponding map $\Wr_\bla: X_\bla \to \C^{(n)}$ is called the Wronski
map.

Denote
\[
X_\bla^0 = X_\bla \cap \Wr_\bla^{-1}\bigl(\big(\C^n-\Delta\big)/S_n\bigr) .
\]

Irreducible representations of the symmetric group $S_n$ are labeled by
partitions $\bla$ of $n$. Denote by $d_\bla$ the dimension of the irreducible
representation corresponding to~$\bla$.
The Wronski map is a f\/inite map of degree~$d_\bla$, see, for example,~\cite{MTV3}.

\subsection[Universal differential operator on $X_{\lambda}$]{Universal dif\/ferential operator on $\boldsymbol{X_{\lambda}}$}

Given an $n\times n$ matrix $A$ with possibly noncommuting entries~$a_{ij}$,
we def\/ine the {\it row determinant\/} to be
\begin{gather*}
\rdet A =
\sum_{\si\in S_n} (-1)^\si a_{1\si(1)}a_{2\si(2)}\cdots a_{n\si(n)} .
\end{gather*}
Let $x=(f_1\lc f_n)$ be a point of $X_\bla$.
Def\/ine the dif\/ferential operator $\D_{\bla,x}$ by
\begin{gather*}
\D_{\bla,x}= \prod_{1\le i<j\le n} \big(\Tilde\la_j-\Tilde\la_i\big)^{-1}\cdot \rdet
 \begin{pmatrix} f_1(u) & f_1'(u) &\dots & f_1^{(n)}(u) \\
f_2(u) & f_2'(u) &\dots & f_2^{(n)}(u) \\ \dots & \dots &\dots & \dots \\
1 & \der &\dots & \der^n
\end{pmatrix} ,
\end{gather*}
where $\der=d/du$. It is a dif\/ferential operator in variable $u$,
\begin{gather}
\label{co diff}
\D_{\bla,x} = \sum_{0\leq i\leq j\leq n} P_{ij}(x)\,u^{n-j}\>\der^{n-i},
\end{gather}
with $P_{ij}\in\C[X_\bla]$.
By formulae (2.11) and (2.3) in~\cite{MTV3},
we have
\begin{gather}
\label{f-la}
\sum_{i=1}^nP_{ii}\prod_{j=i+1}^n(s+j) = \prod_{j=1}^n(s-\la_j+j) ,
\end{gather}
where $s$ is an independent formal variable.

Let  $x\in X_\bla^0$.
 Fix  $\bs z_x=(z_{1,x},\dots,z_{n,x})\in\C^n$ corresponding to
 $\Wr_\bla(x)\in\C^{(n)}$.  Then
\begin{gather*}
\D_{\bla,x} = \prod_{a=1}^n (u-z_{a,x})
\left(
\der^n -\sum_{a=1}^n\frac 1{u-z_{a,x}} \der^{n-1} \right.\\
\left.
\hphantom{\D_{\bla,x} =}{}
+ \sum_{a=1}^n
\frac1{u-z_{a,x}} \left(- p_{a,x}+
\sum_{b\neq a}\frac 1{z_{a,x}-z_{b,x}} \right) \der^{n-2} +\cdots \right)
\end{gather*}
for suitable numbers $\bs p_x=(p_{1,x}\lc p_{n,x})$, see Lemma~3.1 in~\cite{MTV5}.

\begin{lem}\label{X0T}
The map
\[
\psi_\bla: \ X_\bla^0 \to T^*\big(\C^n-\Delta\big)/S_n, \qquad
x \mapsto\bigl(\bs z_x, \bs p_x\bigr) ,
\]
is an embedding whose image is $\Spec/S_n$.
\end{lem}

Lemma~\ref{X0T} is proved in Section~\ref{sec Proof of Lemma }.

\subsection[Description of $\pi^0$]{Description of $\boldsymbol{\pi^0}$}
\label{sec description}

\begin{thm}\label{thm initial}\qquad \null
\begin{enumerate}\itemsep=0pt
\item[$(i)$] The irreducible components of the subvariety $C^0_n\subset C_n$ are naturally labeled
by partitions~$\bla$ of~$n$,  $C^0_n=\cup_\bla C^0_\bla$, where $C^0_\bla$ is the closure of $\xi(\Spec)$ in $C_n$.

\item[$(ii)$] For any $\bla$,
the equations $Q_a=0$, $a=1,\dots,n$, define
$C^0_\bla$ in $C_n$ with multiplicity $d_\bla$.

\item[$(iii)$] The irreducible components of $C^0_n$ do not intersect.
Each component is an $n$-di\-men\-sio\-nal submanifold of $C_n$ isomorphic
to an $n$-dimensional affine space,~{\rm \cite{Wi}}.

\item[$(iv)$] Let $\bla$ be a partition of $n$. Then
there is an embedding $\phi_\bla : X_\bla \to C_n^0$ whose image is
$C^0_\bla$ and such that the following diagram is commutative:
\begin{center}
\begin{picture}(50,50)(0,0)
\put(0,50){\makebox(0,0){$X_\bla$}}
\put(50,50){\makebox(0,0){$C^0_\bla$}}
\put(15,52){\vector(1,0){20}} \put(25,55){\makebox(0,0)[b]{$\phi_\bla$}}
\put(25,0){\makebox(0,0){$\C^{(n)}$}}
\put(7,38){\vector(1,-2){13}}
\put(-4,20){\makebox(0,0)[b]{$\Wr_\bla$}}
\put(43,38){\vector(-1,-2){13}}
\put(45,20){\makebox(0,0)[b]{$\pi^0$}}
\end{picture}
\end{center}
\end{enumerate}
\end{thm}

C.f.~the statements $(i)$--$(iii)$ with results in \cite{FG}.

The map $\phi_\bla$ is given by the following construction. The restriction of
$\phi_\bla$ to $X_\bla^0$ is the composition $\xi\circ\psi_\bla$, where $\xi$
is given by~\eqref{xi}. This map extends from $X_\bla^0$ to an embedding
$X_\bla \to C^0_n$, see Section~\ref{sec end of proofs}.

Parts $(i)$ and~$(ii)$ of Theorem~\ref{thm initial} are proved in Section~\ref{sec part ii proof}.
Parts $(iii)$ and $(iv)$ of Theorem~\ref{thm initial} are proved in Section~\ref{sec end of proofs}.

\begin{rem}
It follows from Theorem~\ref{thm initial} that $\pi^{-1}(0\times 0)$ consists of points
labeled by partitions and the multiplicity of the point corresponding to a partition $\bla$ equals $(d_\bla)^2$.

The fact that the points of $\pi^{-1}(0\times 0)$ are labeled by partitions was explained in~\cite{EG}. The fact that the multiplicity equals $(d_\bla)^2$ was formulated
in~\cite{EG} as Conjecture~17.14 and proved in~\cite{FG}.
\end{rem}

\section{Proofs}
\label{sec Proofs}

\subsection{Proof of Theorem~\ref{thm gener fun}}
\label{Proof of Theorem gener fun}

Assume that a point $\zz=(z_1,\dots,z_n)$ has distinct coordinates.
The Bethe ansatz construction assigns an eigenvector $\omega(\zz,\bs t)$ of
Gaudin Hamiltonians $H_a(\zz)$ on $\sing V^{\otimes n}[\bla]$ to a critical
point $(\zz,\bs t)$ of the master function $\Phi_{\bla}(\zz,\bs t)$,
see \cite{Ba,Ju,MV2,RV},
\begin{gather}
\label{eig}
H_a(\zz) \omega(\zz,\bs t) = \frac {\der \Phi_{\bla}}{\der z_a}
(\zz, \bs t) \omega(\zz, \bs t) ,\qquad a=1,\dots,n.
\end{gather}
Formula \eqref{eig} shows that $L_\bla\subset\Spec$.
By Theorem 6.1 in \cite{MV2} the Bethe vectors form a basis of
$\sing V^{\otimes n}[\bla]$ for generic $\zz \in \C^n-\Delta$. This proves Theorem~\ref{thm gener fun}.

\subsection[Proof of part (i) of Theorem \ref{thm spec-lambda in L_0}]{Proof of part $\boldsymbol{(i)}$ of Theorem \ref{thm spec-lambda in L_0}}
\label{Proof of part i}
It is easy to see
that $\sing V^{\otimes n}[\bla]$ and the action on it of the
Hamiltonians $H_a(\zz)$ do not depend on $N$. Hence, $\Spect_{N,\bla}$
does not depend on $N$.

Another (less straightforward) proof of part $(i)$ follows from
formula \eqref{eig} and the fact that~$\Phi_{N,\bla}$ does not depend on~$N$.

\subsection{Proof of Lemma~\ref{X0T}}
\label{sec Proof of Lemma }

\begin{lem}
\label{lem spec smooth}
For every $\bla$, the spectral variety $\Spec/S_n \subset T^*(\C^n-\Delta)/S_n$ is
smooth. For different $\bla$'s the spectral subvarieties do not intersect.
\end{lem}

\begin{proof}
Let  $x\in X_\bla^0$ and  $\Wr_\bla(x)$
be a projection of $\zz_x=(z_{1,x},\dots,z_{n,x})$.
By~\cite{MTV3}, the points $x\in X_\bla^0$ are in
a one-to-one correspondence with the eigenvectors of the Gaudin
Hamiltonians on $\sing V^{\otimes n}[\bla]$. Denote $v_x$ the eigenvector
corresponding to $x$. By Lemma~3.1 in~\cite{MTV5} the numbers
$\bs p_x=(p_{1,x},\dots,p_{n,x})$ are eigenvalues of
$H_1(\zz_x)\lc H_n(\zz_x)$ on $v_x$. Hence, the image of $\psi_\bla$ is
$\Spec/S_n$.

By Theorem~3.2 in~\cite{MTV5}, the coordinates $z_{1,x}\lc z_{n,x}$,
$p_{1,x}\lc p_{n,x}$ generate all functions on~$X_\bla^0$. This proves that
$\psi_\bla$ is an embedding of $X_\bla^0$ to $T^*(\C^n-\Delta)/S_n$.

By Theorem 3.2 in~\cite{MTV5}, the coef\/f\/icients $P_{ij}(x)$ in~\eqref{co diff}
are given by some universal functions in $\zz_x$, $\bs p_x$ independent of~$\bla$. Hence formula~\eqref{f-la} implies that the spectral varieties for
dif\/ferent $\bla$'s do not intersect.
\end{proof}

\begin{lem}
\label{lem spec in L0}
For every $\bla$ the spectral variety~$\Spec$ lies in the variety~$L_0$ defined
in~\eqref{L0}.
\end{lem}

\begin{proof}
Let  $x\in X_\bla^0$ and $\bs z_x=(z_{1,x},\dots,z_{n,x})$ corresponds to  $\Wr_\bla(x)$.
By Theorem 3.2 in~\cite{MTV5},
\begin{gather*}
\det\bigl((u-Z_x)(v-Q_x)-1\bigr) =
\sum_{0\leq i\leq j\leq n} P_{ij}(x) u^{n-j} v^{n-i} ,
\end{gather*}
where $Z_x = \diag(z_{1,x},\dots,z_{n,x})$, $Q_x$ is given by \eqref{Q} in terms of~$\bs z_x$ and~$\bs p_x$, and
$P_{ij}(x)$ are given by \eqref{co diff}. This implies that $\det(v-Q_x)=v^n$
and hence~$\Spec\subset L_0$.
\end{proof}

\subsection[Proofs of parts (i) and (ii) of Theorem~\ref{thm initial},
part (ii) of Theorem~\ref{thm spec-lambda in L_0}, and
Theorem~\ref{thm L_0 smooth}]{Proofs of parts $\boldsymbol{(i)}$ and $\boldsymbol{(ii)}$ of Theorem~\ref{thm initial},\\
part~$\boldsymbol{(ii)}$ of Theorem~\ref{thm spec-lambda in L_0}, and
Theorem~\ref{thm L_0 smooth}}
\label{sec part ii proof}

Consider the Lie algebra $\frak{gl}_n$ with standard generators $e_{ij}$.
Fix a set of complex numbers $\bs q=(q_1,\dots,q_n)$. Consider the
weight subspace
\[
V^{\otimes n}[1,\dots,1]=\big\{w\in V^{\otimes n}\; | \; e_{ii}w=w,\ i=1,\dots,n\big\}
\]
and the Gaudin Hamiltonians
\begin{gather*}
H_a(\bs z,\bs q)   =  \sum_{i=1}^n q_i e_{ii}^{(a)} +
\sum_{i,j=1}^n \sum_{b\neq a}  \frac{e_{ij}^{(a)}e_{ji}^{(b)}}{z_a-z_b} ,
\end{gather*}
which generalize the Gaudin Hamiltonians in~\eqref{Ha}.
The generalized Gaudin Hamiltonians act on $V^{\otimes n}[1,\dots,1]$.
By Theorem~5.3 in~\cite{MTV1}, for generic $(\bs z,\bs q)$ the
generalized Gaudin Hamiltonians
$H_1(\bs z,\bs q), \dots, H_n(\bs z,\bs q)$ have an eigenbasis in $V^{\otimes n}[1,\dots,1]$. By Theorem~4.3
in~\cite{MTV4} the eigenvectors are in one-to-one correspondence with
the preimages of the point $(\zz,\bs q)$ under the map
$\pi : C_n \to \C^{(n)}\times\C^{(n)}$.
This identif\/ication sends an eigenvector $w$
with  $H_a(\bs z,\bs q)w=p_aw$, $a=1,\dots,n$, to the point $(Z,Q)$,
where $Z=\text{diag}(z_1,\dots,z_n)$ and $Q$ is def\/ined by~\eqref{Q}.

The $\frak{gl}_n$-action on every eigenvector of the Gaudin Hamiltonians
$H_a(\bs z,\bs q=0)$ on the space $\sing V^{\otimes n}[\bla]$ generates a
$d_\bla$-dimensional subspace in $V^{\otimes n}[1,\dots,1]$ of eigenvectors of
\mbox{$H_a(\bs z,\bs q{=}0)$}. The identif\/ication
of Theorem~4.3 in~\cite{MTV4} implies that $\Spec$ has multiplicity $d_\bla$
when def\/ined by the
equations $Q_a(\zz,\bs p)=0$, $a=1,\dots, n$.
More precisely, the identif\/ication tells that $d_\bla$ points of the $\pi^{-1}(\bs z,\bs q)$
collide to one point of $\Spec$, when $\bs q\to 0$.
This proves part $(ii)$ of Theorem~\ref{thm initial}.

Since $\sum_\bla d^2_\bla =n!$, we conclude that $L_0=\cup_\bla\Spec$. This proves part~$(i)$ of
Theorem \ref{thm initial}, part~$(ii)$ of Theorem~\ref{thm spec-lambda in L_0},
and Theorem \ref{thm L_0 smooth}.

\subsection[Proof of parts (iii) and (iv) of Theorem \ref{thm initial}]{Proof of parts $\boldsymbol{(iii)}$ and $\boldsymbol{(iv)}$ of Theorem \ref{thm initial}}
\label{sec end of proofs}

Every point of $X_\bla$ is a point of Wilson's adelic Grassmannian $\text{Gr}^{\text{ad}}(n)$,
which is identif\/ied with the Calogero--Moser space $C_n$ by the main Theorem~5.1
in~\cite{Wi}. For generic points of $X_\bla$ the map is $x\mapsto (Z_x,Q_x)$.
The induced map of functions $\C[C_n]\to \C[X_\bla]$ is constructed as follows.
Consider
$P(u,v)=\det((u-Z)(v-Q)-1)$
as a polynomial in $u,v$ whose coef\/f\/icients are functions on~$C_n$,
$P(u,v)=\sum_{i,j}\Tilde P_{ij}u^{n-j}v^{n-i}$.
The coef\/f\/icients $\Tilde P_{ij}$ generate $\C[C_n]$, see Lemma~4.1
in~\cite{MTV4}.
The map $\C[C_n]\to \C[X_\bla]$ is def\/ined by the formula
$\Tilde P_{ij} \mapsto P_{ij}$,
where  $P_{ij}$ are given by~\eqref{co diff}, see Theorem 4.3 in~\cite{MTV4}.
By Lemma~3.4 in~\cite{MTV3} the image of this map is $\C[X_\bla]$. Hence
$X_\bla \to C_n$ is an embedding. By formula~\eqref{f-la} the images
do not intersect for dif\/ferent~$\bla$'s. This proves parts~$(iii)$ and~$(iv)$
of Theorem~\ref{thm initial}.

\section{Further remarks}
\label{sec Farther remarks}

Fix distinct complex numbers $\bs q=(q_1,\dots,q_n)$. Let
$\sigma_a(\bs q)$, $a=1,\dots,n$, be the elementary symmetric functions of
$\bs q$. Def\/ine
\begin{gather*}
L_{\bs q} = \big\{(\bs z,\bs p)\in T^*(\C^n-\Delta)\; |\; Q_a(\bs z, \bs p)=\sigma_a(\bs q), \ a=1,\dots,n\big\} .
\end{gather*}

\begin{thm}
\label{thm L_q smooth}
The subvariety $L_{\bs q}$ is irreducible, smooth, and Lagrangian.
\end{thm}

Def\/ine the spectral variety of the Gaudin Hamiltonians $H_a(\zz, \bs q)$, $a{=}1,\dots,n$, on
$V^{\otimes n}[1,\dots,1]$ by the formula
\begin{gather*}
\Spect_{\bs q} = \big\{(\bs z,\bs p)\in T^*(\C^n-\Delta)\; |\; \exists\,
v\in V^{\otimes n}[1,\dots,1]
\ \text{with}\ H_a(\bs z, \bs q)v = p_av,\ a=1,\dots,n\big\}.
\end{gather*}

\begin{thm}
\label{thm spec-lambda in L_q}
We have $L_{\bs q} = \Spect_{\bs q}$.
\end{thm}

Consider the set of $n(n-1)/2$ variables
\[
\bs t = \big(t^{(1)}_{1},\dots,t_{n-1}^{(1)},\dots, t^{(n-2)}_{1},
t^{(n-2)}_{2},t^{(n-1)}_{1}\big)
\]
and the af\/f\/ine space $\C^n\times\C^{n(n-1)/2}$ with coordinates
$\bs z$, $\bs t$. Consider the master
function $\Phi_{\bs q} : \C^n\times\C^{n(n-1)/2}\to\C$,
\begin{gather*}
\Phi_{\bs q} (\bs z,\bs t) =
\sum_{1\leq a<b\leq n}  \log(z_a-z_b)
-\sum_{a=1}^{n}\sum_{i=1}^{n-1} \log\big(t_i^{(1)}-z_a\big)
\\
\hphantom{\Phi_{\bs q} (\bs z,\bs t) = }{}
+ 2 \sum_{k=1}^{n-1}\sum_{1\leq i<j\leq n-k} \log \big(t_i^{(k)}-t_j^{(k)}\big)-
\sum_{k=0}^{n-2} \sum_{i=1}^{n-k} \sum_{j=1}^{n-k-1}
\log\big(t_i^{(k)}-t_j^{(k+1)}\big)
\\
\hphantom{\Phi_{\bs q} (\bs z,\bs t) = }{}
+\sum_{k=1}^{n-1}\sum_{i=1}^{n-k} (q_{k+1}-q_k)t^{(k)}_i
+ q_1\sum_{a=1}^n z_a ,
\end{gather*}
see~\cite{FMTV}.
Denote by $\Crit_{\bs q}$ the critical set of $\Phi_{\bs q}$ with respect to
$\bs t$,
\[
\Crit_{\bs q} = \big\{(\bs z,\bs t)\in \C^n\times\C^{n(n-1)/2}\ | \ d_\TT\Phi_{\bs q}(\bs z,\bs t) = 0\big\} .
\]
Denote by $\Tilde L_{\bs q}\subset T^*(\C^n-\Delta)$ the image of the map
\begin{gather*}
\Crit_{\bs q} \to T^*(\C^n-\Delta), \qquad (\bs z,\bs t) \mapsto (\bs z,\bs p), \qquad
\text{where}\quad p_a = \frac {\der \Phi_{\bs q}}{\der z_a}(\bs z,\bs t),\quad a=1,\dots,n.
\end{gather*}

\begin{thm}
\label{thm-gener-fun}
We have $\Tilde L_{\bs q}\subset\Spect_{\bs q}$ and
$\Spect_{\bs q}$ is the closure of $\Tilde L_{\bs q}$ in $T^*(\C^n-\Delta)$.
\end{thm}

Consider the canonical map $ \pi : C_n \to \C^{(n)}\times \C^{(n)},$ $ (Z,Q)
\mapsto (\text{spec}(Z), \text{spec}(Q))$.
Denote by $C_n^{\bs q}$ the subvariety $\pi^{-1}(\C^{(n)}\times \bs q)
\subset C_n$. Identifying $\C^{(n)}\times \bs q$ with $\C^{(n)}$ we get a map
$ \pi^{\bs q} : C^{\bs q}_n \to \C^{(n)},$ $(Z,Q) \mapsto \text{spec}(Z)$,
induced by~$\pi$.

Denote
\[
f_i(u) = e^{q_iu}(u+f_{i1}), \qquad
i=1,\dots,n.
\]
Denote by $X_{\bs q}$ the $n$-dimensional af\/f\/ine space of $n$-tuples $\{f_1,\dots,f_n\}$
of such quasiexponentials.
The polynomial algebra
$\C[X_\bla]\,=\,\C[f_{11},\dots,f_{n1}]$
is the algebra of regular functions on $X_{\bs q}$.
We have
\begin{gather*}
\Wr(f_1(u),\dots,f_n(u)) = e^{(q_1+\dots+q_n)u}
\prod_{1\le i<j\le n} (q_j-q_i)\cdot
\left(u^n+\sum_{a=1}^n (-1)^a W_a u^{n-a}\right)
\end{gather*}
with $W_1,\dots,W_n\in\C[X_{\bs q}]$.
Def\/ine an algebra homomorphism
\[
\W_{\bs q} : \ \C[\C^{(n)}] \to \C[X_{\bs q}] ,\qquad
\sigma_a \mapsto W_a .
\]
Let $\Wr_{\bs q} : X_{\bs q} \to \C^{(n)}$ be the corresponding map of spaces.

\begin{thm}\label{thm initial q} \qquad \null
\begin{enumerate}\itemsep=0pt
\item[$(i)$] The equations $Q_a=q_a$,  $a=1,\dots,n$, define
$C^{\bs q}_n$ in $C_n$ with multiplicity~$1$.

\item[$(ii)$] There is an embedding $\phi_{\bs q} : X_{\bs q} \to C_n^{\bs q}$
whose image is $C^{\bs q}_n$ and such that the following diagram is
commutative:

\begin{center}
\begin{picture}(50,50)(0,0)
\put(0,50){\makebox(0,0){$X_{\bs q}$}}
\put(50,50){\makebox(0,0){$C^{\bs q}_n$ }}
\put(15,52){\vector(1,0){20}} \put(25,55){\makebox(0,0)[b]{$\phi_{\bs q}$}}
\put(25,0){\makebox(0,0){$\C^{(n)}$}}
\put(7,38){\vector(1,-2){13}}
\put(-4,20){\makebox(0,0)[b]{$\Wr_{\bs q}$}}
\put(43,38){\vector(-1,-2){13}}
\put(45,20){\makebox(0,0)[b]{$\pi^0$}}
\end{picture}
\end{center}
\end{enumerate}
\end{thm}

The map $\phi_{\bs q}$ is given by the following construction.
For  $x=(f_1\lc f_n)\in X_{\bs q}$, def\/ine the dif\/ferential operator
$\D_{\bs q,x}$ by
\begin{gather*}
\D_{\bs q,x}= e^{-(q_1+\dots+q_n)u}\prod_{1\le i<j\le n} (q_j-q_i)^{-1}
\rdet  \begin{pmatrix} f_1(u) & f_1'(u) &\dots & f_1^{(n)}(u) \\
f_2(u) & f_2'(u) &\dots & f_2^{(n)}(u) \\ \dots & \dots &\dots & \dots \\
1 & \der &\dots & \der^n
\end{pmatrix} .
\end{gather*}
Then
\[
\D_{\bs q,x}= \sum_{i=0}^n \sum_{j=0}^n
P_{ij}(x) u^{n-j} \der^{n-i} ,
\]
where $P_{ij}\in\C[X_{\bs q}]$.

Denote
$X_{\bs q}^0= X_{\bs q}\cap \Wr_\bla^{-1}\bigl((\C^n-\Delta)/S_n\bigr)$
 and consider the map
\[
\psi_{\bs q}: \ X_{\bs q}^0 \to T^*\big(\C^n-\Delta\big)/S_n, \qquad
x \mapsto\bigl(\bs z_x , \bs p_x\bigr) ,
\]
where $\bs z_x\in\C^n$ projects to  $\Wr_{\bs q}(x)\in\C^{(n)}$ and
 $\bs p_x=(p_{1,x}\lc p_{n,x})$,
\[
p_{a,x} = -\Res_{u=z_{a,x}}
\left( \frac{\sum\limits_{j=0}^n P_{2,j}(x) u^{n-j}}
{\prod\limits_{i=1}^n (u-q_i)}\right) + \sum_{b\neq a}\frac 1{z_{a,x}-z_{b,x}} .
\]
Then the restriction of $\phi_{\bs q}$ to $X_{\bs q}^0$ is the composition
$\xi\circ\psi_{\bs q}$, where $\xi$ is given by~\eqref{xi}. This map extends
from $X_{\bs q}^0$ to an embedding $X_{\bs q} \to C^{\bs q}_n$.

The proofs of Theorems~\ref{thm L_q smooth}--\ref{thm initial q}
are basically the same as the proofs of Theorems~\ref{thm L_0 smooth}--\ref{thm initial}.

\subsection*{Acknowledgments}

The authors thank V.~Schechtman and D.~Arinkin for useful discussions.
The third author thanks the MPIM, HIM, and IHES for the hospitality.
E.~Mukhin was supported in part by NSF grant DMS-0900984.
V.~Tarasov was supported in part by NSF grant DMS-0901616.
A.~Varchenko was supported in part by NSF grants DMS-0555327 and DMS-1101508.

\pdfbookmark[1]{References}{ref}
\LastPageEnding

\end{document}